\documentclass[]{article}
\usepackage{amsfonts}
\title{Growth rates of small
cancellation groups.}
\author{Anna Erschler (Dyubina), Tel Aviv University}
\date{e-mail: annadi@tau.ac.il,
 erschler@pdmi.ras.ru}

\begin{document}
\newcommand{\R}{{\mathbb R}}
\newcommand{\N}{{\mathbb N}}
\newcommand{\Z}{\mathbb Z}
\newcommand{\con}{{\rm Con}}
\newcommand{\rad}{{\rm R}}
\newcommand{\dist}{{\rm dist}}
\newcommand{\cyl}{{\rm Cyl}}
\renewcommand{\phi}{\varphi}
\newcommand{\ti}{\tilde}
\newtheorem{lemma}{Lemma}
\newtheorem{definition}{Definition}
\newtheorem{example}{Example}
\newtheorem{cor}{Corollary from corollary}
\newtheorem{corollary}{Corollary}
\newtheorem{theorem}{Theorem}
\newtheorem{notation}{Notation}
\newtheorem{proposition}{Proposition}
\newtheorem{rem}{Remark}
\maketitle

\begin{abstract}
We estimate growth rate for small cancellation groups.
In particular we show that there is a continuum of possible
values for exponential growth rates.
\end{abstract}

\section{Introduction}

Let $G$ be a group with some fixed set of generators and let $d$ denote
the word metric corresponding to this set of generators.

Let $B_G(n)$ denote the ball $\{g\in G: d(g,e)\le n\}$
and $v(G)$ denote the (exponential) growth rate
$$
v(G)=\lim\sqrt[n]{\#B_G(n)}.
$$
The group $G$ has exponential growth, if for some (and hence for all)
set of generators of $G$ $v(G)>1$. There are many examples
of groups of exponential growth (e.g. free groups and more generally
any hyperbolic, any non-amenable group; any solvable but not
virtually nilpotent group \cite{Mil}, \cite{Wolf}
). Moreover, all known examples of
finitely presented groups of non-exponential growth have polynomial
growth
(and hence are virtually nilpotent by a theorem of
Gromov \cite{Gromov2}).
We remind that there exist groups of
intermediate growth, but all known examples
are not finitely presented. The existence  of such groups
was first discovered by R.Grigorchuk
 \cite{Gri}.

Yet not much  is known about values of exponential growth rates.
Set
$$
\Omega^G=\{v\in [1,\infty) | v=v(G)
\mbox{ for  some set of generators  of } G \}.
$$
A well known question \cite{GrPa}
which remains open for more than 20 years asks if there exists a group $G$
of exponential growth such
that
$$
\inf v_{v \in \Omega^G}  =1.
$$
That is if there exists a group of exponential, but not uniformly
exponential growth.
Recently some partial results for this question
were obtained (see references in \cite{Ha}). Namely it was proved 
for some classes of groups that
exponential growth implies uniformly exponential
growth.

Another series of questions about growth rates is which number can appear
as the growth rate of some group (or for some group in the given class of
groups).

If we fix the number $n$ of generators of $G$ then there are obvious
estimates
$$
1 \le v \le 2n-1,
$$
where the equality in the last inequality is achieved only if
the group is free on the $2n-1$ given generators.
There are no other known   general restrictions for $v$.
 But for some classes of groups there exist such restrictions.
  For example
any growth rate of any hyperbolic group is an interger algebraic 
number \cite{Gromov},
\cite{GyHa}.

The groups we consider in this paper are close to hyperbolic ones.
We consider not finitely presented groups with small cancellation
law. We remind that having a finite presentation would imply hyperbolicity.
Still the situation for these groups differs from the hyperbolic case.

In this paper we study the set
$$
\Omega=\{v\in [1,\infty) | v=v(G) \mbox{ for some  }G \}.
$$
It is known (\cite{GriHa}, \cite{Ha}) that $\Omega$ is everywhere
dence in $[1,\infty]$.
It was asked in \cite{Ha} whether $\Omega$ has the cardinality
of continuum.

The following theorem answers this question.
\begin{theorem}
Let $\Omega_2=\{v\in [1,\infty) | v=v(G)$ for some
two generated group $G \}$. Then

 $\# \Omega_2 =2^{\aleph_0}$.
\end{theorem}

This theorem follows from the following proposition

\begin{proposition}
There exists $E:\N \to \N$ with the following property.
Let $r_i=(a^{E(i)}b^{E(i)})^{100}$.
For any $J\subset \N$ define $G_J$ by
$$
G_J=\langle a,b|r_j=e, j\in J \rangle   .
$$

Then $I\ne J$ implies that $v(G_{I})\ne v(G_{J})$.

\end{proposition}

\begin{rem}
Similar groups were first introduced by B.H.Bowditch. He used them to construct
continuously many quasi-isometry classes of groups \cite{Bow}.
Similar groups provide also examples of a group with non-homeomorphic
asymptotic cones \cite{TV}.
\end{rem}

\section{Definitions and preliminary observations}

Let $w$ be a word and $|w|$ denote its length in the corresponding
free group.

For a non-negative matrix $A$
it is known (\cite{Gant}) that there exists a non-negative eigenvalue $v$
of $A$ such that for any other eigenvalue $v'$ $\|v'\|\le v$.
Let
 $v(A)$ denote this maximal
eigenvalue and $\sum (A)$ denote the sum of elements of $A$.

Let $G$ be a group and let $\langle S,R \rangle$
 be a presentation of $G$.
{\it Symmetrization} $R_*$ consists of all distinct cyclic permutations
of the defining relators $r$ and their inverses.
A word $u$ is a {\it piece} relative to $R$, if $R_*$ contains two distinct
elements of the form $uv'$ and $uv''$.

\begin{definition} (See, for example, \cite{Stre})
Let $\lambda>0$.
We say that $S,R$ satisfy {\it condition} $C'(\lambda)$ if
for any $r\in R_*$ and any piece $u$ of it
$$|u|<\lambda |r|.$$
\end{definition}

If some finite presentation of a group satisfies $C'(1/6)$, then
this group is hyperbolic \cite{Stre}. Hence there is a Markov grammar
corresponding  to $G$ (and preserving the word length) \cite{Gromov},
\cite{GyHa}.

\begin{lemma}
 Let
$$
H= \langle S|w_i=e, i \in \N \rangle
$$
and
$$
H_k= \langle S| w_i=e, i\le k \rangle.
$$
Then $\lim v(H_k)=v(H)$.
\end{lemma}
{\bf Proof.}
It is clear that $v(H_k)\ge v(H_j)$ for $k< j$.
Hence there exists a limit $v=\lim v(H_k)$.
Then $v(H)\le v$, since $v(H)\le v(H_k)$ for any $k$.

On the other hand
$$
\#B_{H_k}(r)\ge v(H_k)^r\ge v^r
$$
and for each $r$ there exists $k$ such that
$$
\# B_H(r)=\# B_{H_k}.
$$
Hence for any $r$ $\# B_H(r)\ge v^r$ and $v(H)\ge v$.

\begin{definition}
Consider a Markov grammar with states $1,2,...,N$.
Adjacency matrix of this grammar $M$ is the matrix such that
$M(i,j)$ is the number of oriented edges from state $i$ to state $j$.
\end{definition}

If we say that a Markov grammar corresponds to some group
we assume that the correspondence preserves word length.

\begin{lemma} Consider a language given by some Markov grammar
and let $M$ be the adjacency matrix of this grammar.
Then growth rate of the language is equal to $v(M)$.
\end{lemma}

{\bf Proof.}
It is clear that growth rate is not greater than $v(M)$.
On the other hand we know that there is a non-negative eigenvector
corresponding to $v(M)$ (\cite{Gant}). This implies that
$\sum (M^r) \ge Cv(M)^r$ and so the growth rate is not less than $v(M)$.

\section{Estimate from below}

\begin{lemma} Consider a Markov grammar with $s$ states and $k$ labels.
Let $v$ be the growth rate of the corresponding language.
Suppose that $s\le N$ and fix some $p$ words of length $4N$.
We call them long forbidden words.
Consider the following new language. It consists of  the words of the old
language that do not contain any of the long forbidden words.
Let $v_{new}$ be the growth rate of the new language.
Then
$$
v_{new}^N\ge (v^N/s^2-4Np)/k^s
$$
\end{lemma}

{\bf Proof.}
Let $M$ be the adjacency matrix of the Markov grammar.
Then $p_{i,j}^L=M^L(i,j)$ is the number of paths of length $L$ starting
from the state $i$ and ending in the state $j$.

Note that there exists $\tilde{s}\le s$ such that for some
$i$

$$
p_{i,i}^{N+\tilde{s}} \ge v^N/s^2.
$$

In fact, let $A_1, A_2,... A_f$ be irreducible blocks of $M$.
Then there exists $g$ such that $v(A_g)=v(M)=v$.
Since $\sum(A_g^N)\ge v(A_g)^N=v^N$ we see that
there exist $i$ and $j$ in the block $A_g$ such
that
$$
p_{i,j}^N=M^N(i,j)=A_g^N(i,j)\ge v^N/s^2.
$$

There exists $\tilde{s}\le s$ such that we can go from the state $i$
to the state $j$ by $\tilde{s}$ steps.
This implies that
$$
p_{i,i}^{N+\tilde{s}} \ge v^N/s^2.
$$

Let $N'=N+\tilde{s}$. Note that $N'\le 2N$.
Consider all subwords of length $N'$ of long forbidden words.
We will call them short forbidden words.
Note that there are at most $4Np$ short forbidden words.
Take some word $w$. Suppose that
$$
w=a_1a_2...a_l
$$
and that $|a_m|=N'$ for any $m<l$, $|a_l|\le N'$.
Note that if none of $a_m$ is a short forbidden word, then
$w$ does not contain any of long forbidden words.

This implies that
$$
v_{new}^{N'}\ge v^N/s^2-4Np.
$$

Since $v_{new}\le v\le k$ we deduce that
$$
v_{new}^N\ge (v^N/s^2-4Np)/k^s.
$$
\begin{lemma}
Let $$
G=\langle S|w_1=w_2=...=w_{n+1}=e \rangle
$$
and
$$
H=\langle S|w_1=w_2=...=w_n=e \rangle.
$$
Suppose that these presentations satisfy $C'(1/6)$ and that
$|w_i|>6$ for any $i$.
We call a word $w$ $1/12$-reduced with respect to $w_{n+1}$ if it contains
no subword of length $\ge 1/12|w_{n+1}|$ of $w_{n+1}$ and if
it is freely reduced.
\begin{enumerate}
\item If $w$ be $1/12$-reduced with respect to $w_{n+1}$ and
$w$ is geodesic in $H$,  then it is geodesic in $G$.

\item Let $\gamma_1$ and $\gamma_2$ be $1/12$-reduced with respect to
$w_{n+1}$. If $\gamma_1$ and $\gamma_2$ are geodesic in $H$ and
$\gamma_1\ne \gamma_2$ in $H$ then $\gamma_1\ne \gamma_2$ in $G$.
\end{enumerate}
\end{lemma}
{\bf Proof.}
\begin{enumerate}
\item Let us use induction on $|w|$.

{\bf Base.} Note that if $|w|=1$, then the $C'(1/6)$ condition
implies that
$w\ne e$ in $G$ and hence
it is geodesic in $G$.

{\bf Induction step.}
Suppose that $w=w_1w_2$ and that $|w_2|=1$.
Since $|w_1|=|w|-1$ and $w_1$ is geodesic in $H$, we know already
that $w_1$ is geodesic in $G$.
Suppose that $w=w_3$ in $G$ and that $|w_3|<|w|$.
Consider a triangle in $G$ with sides $w_1, w_2, w_3$.
Without loss of generality we can assume that this is a simple
reduced triangle, that is that its sides intersect only in its vertices.
In fact, otherwise we can replace $w_1$ and $w_3$ with some of its
end subwords.
Consider a reduced diagram of this triangle and its tiling.
Note that the tiling contains at most 2 distinguished tiles.
(A tile is said to be distinguished, if it contains a vertex of
a triangle inside some of its open exterior edges (see \cite{Stre})).
Then part 3 of theorem 35 of \cite{Stre} implies that
the tiling looks like one on the picture.

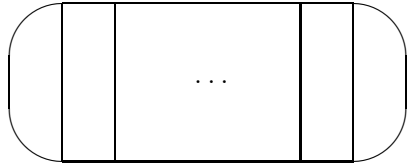
\begin{picture}(300,80)
\put(150,40){\oval(150,60)}
\put(95,10){\line(0,1){60}}
\put(115,10){\line(0,1){60}}
\put(205,10){\line(0,1){60}}
\put(185,10){\line(0,1){60}}
\put(150,40){.}
\put(145,40){.}
\put(155,40){.}
\end{picture}

Note that if some tile corresponds to $w_{n+1}$, then it has
at least 3 interior edges. So no tile can be labelled with $w_{n+1}$.
This means that we have this triangle not only in $G$ but also in $H$.
So $w=w_3$ in $H$, but this is impossible.

\item
By the first part of the lemma we know that $\gamma_1, \gamma_2$ are
geodesic in $G$.
Suppose that $\gamma_1=\gamma_2$ in $G$.
Consider a digon with sides $\gamma_1$ and $\gamma_2$. Again without
loss of generality we can assume that the digon is reduced and
apply part 3 of theorem 35 of \cite{Stre}.
No tile can be labeled with $w_{n+1}$ and this implies
that $\gamma_1=\gamma_2$ in $H$.
\end{enumerate}

Two previous lemmas imply the following corollary.
\begin{corollary}
 Let $$
G=\langle S|w_1=w_2=...=w_{n+1}=e \rangle
$$
and
$$
H=\langle A|w_1=w_2=...=w_n=e \rangle.
$$
Suppose that these presentations satisfy $C'(1/6)$ and
 that there is a Markov grammar  for $H$ that consists of
$s$ states.
There exists $\gamma(s,k)$ such that if
$$
|w_{n+1}|>\gamma(s,k),
$$
then
$$
v(G)\ge v(H)-200/\sqrt{|w_{n+1}|}
$$
\end{corollary}
{\bf Proof.}
Let $N$ be the integer part of $|w_{n+1}|/52$.
We can assume that $N>1000s^2(2k)^{s+6}$.

Consider all subwords of length $4N$ of $w_{n+1}$.
There are at most $52(N+1)$ such subwords.
By the previous lemma we can estimate the growth rate of $G$ by the growth
rate of the grammar corresponding to $H$ where these words of length
$4N$ are forbidden.
Then we apply  lemma 3 and get that
$$
v(G)^N\ge (v(H)^N/s^2-4N 52(N+1))/(2k)^s \ge (v(H)^N-N^4)N.
$$

Note that without loss of generality we can assume that
$v(H)>N^4+1$, since otherwise
 $$v(H)<\sqrt[N]{1+N^4}\le 1+1/\sqrt{N}$$
(note that $N > 1000$).

But if $v(H)>N^4+1$ then
$$\sqrt[N]{v(H)^N-N}\ge v-1/\sqrt{N}$$
and
$$\sqrt[N]{N}\le 1+1/N^{\frac{2}{3}}.$$
Hence
$$
v(G)\ge \frac{v(H)-1/\sqrt{N}}{1+1/N^\frac{2}{3}}\ge
(v(H)-1/\sqrt{N})(v(H)-1/N^\frac{2}{3}).
$$
But since $\sqrt[6]{N}\ge 2k \ge v(H)$
$$
(v(H)-1/\sqrt{N})(v(H)-1/N^\frac{2}{3}
\ge v(H)-2/\sqrt{N}\ge v(H)-200/\sqrt{|w_{n+1}|}
$$

The previous corollary implies another corollary:
\begin{corollary}
Let
 $$
H=\langle S|w_1=w_2=...=w_n=e \rangle
$$
 and
$$
H_k=\langle S|w_1=w_2=...=w_n=\tilde{w}_k=e \rangle.
$$
Suppose that these presentations satisfy $C'(1/6)$ and
 that $|\tilde{w}_k|\to \infty$ when $k\to \infty$.
Then $\lim v(H_k)=v(H)$
\end{corollary}
This generalizes a corollary from \cite{Sh} where a free group
was taken for $H$.

\section{Estimate from above. Decrease of growth rate.}
\begin{notation}
\begin{enumerate}
\item
Choose $\alpha:\N\to \N$ such that
if
$$
H_1=\langle a,b |r_i^1=e, i\in I \rangle,
$$
$$
H_2=\langle a,b |r_j^2=e, j\in J \rangle
$$
and
 $|r_i^1|, |r_j^2|\le N$
then either $v(H_1)=v(H_2)$ or $\|v(H_1)-v(H_2)\| \ge 1/\alpha(N)$.

\item
 Choose $\beta:\N\to \N$ such that
if
$$
H=\langle a,b |r_i=e, i\in I \rangle,
$$
$|r_i|\le N$ for any $i$ and there is a Markov grammar for $H$ (and
for these generators $a,b$), then one of such Markov grammars contains
at most $\beta(N)$ states (note that for example we can take
any $C'(1/6)$ group for $H$).

\end{enumerate}
\end{notation}
Note that $\alpha$ and $\beta$ are well defined, since there are
finitely many groups $H$ such that
$$
H=\langle a,b |r_i=e, i\in I \rangle,
$$
and $|r_i|\le N$ for any $i$.

Just the existence of $\alpha$ will be sufficient for
our estimate from above if we know  that $v(H_1)\ne v(H_2)$.
So the main thing we have to show is that for certain $H_1, H_2$
the growth rates are not equal.
 Namely we will prove the following proposition.

\begin{proposition}
 Let $$
G=\langle S|w_1=w_2=...=w_n=(a^Nb^N)^{100}=e \rangle
$$
and
$$
H=\langle A|w_1=w_2=...=w_n=e \rangle.
$$
Suppose that $N>|w_i|$ for any $i$.

Suppose also that these presentations satisfy $C'(1/50)$ and that none of
$a^{\pm 1}b^{\pm 1}a^{\mp 1}$ and $b^{\pm 1}a^{\pm 1}b^{\mp 1}$
(for any choice of signs)
is a piece for these presentations.
Then $v(G)<v(H)$.
\end{proposition}
\begin{rem} It is conjectured in \cite{GriHa2} that
the growth rate decreases by taking a
quotient of  any hyperbolic
group over any infinite subgroup .
 The proposition 2 would have followed from
this conjecture.
\end{rem}
First we prove the following lemma.

\begin{lemma}
Suppose the assumptions of the proposition hold.
In this lemma we consider the group $H$, so all equalities
and statements about words correspond to $H$.
\begin{enumerate}

\item Let $\gamma$ be a geodesic word  ending with
$a^{\pm 1}b^{\pm 1}$ (for some choice of signs), that is
$\gamma=\gamma'a^{\pm 1}b^{\pm 1}$.
Then
$\gamma(a^N b^N)^k$ is geodesic for any positive $k$.
Moreover, any geodesic from $e$ to $\gamma(a^N b^N)^k$ passes
through $\gamma$ and coincides from this point with
$\gamma(a^N b^N)^k$.

\item
Suppose that $\gamma_1=\gamma_1'b^k$ is geodesic and $k\ge N$.
Let $\gamma_2$ be geodesic that begins with  $a^{\pm 1}b^{\pm 1}$,
that is $\gamma_2=a^{\pm 1}b^{\pm 1}\gamma_2'$.
Then $\gamma_1\gamma_2$ is geodesic. Moreover, any geodesic
from $e$ to $\gamma_1\gamma_2$  passes through $\gamma_1$.

\end{enumerate}
\end{lemma}
{\bf Proof of the lemma.}
Take $\alpha>0$. We say that some word $u$ is $>\alpha W$ if
it is a piece  of some $w_i$ and $|u|/|w_i|>\alpha$.

First note that
from Greendlinger lemma \cite{LynS}
it follows that
 $(a^Nb^N)^k$ and $b^k$ are geodesic for any $k$ and that
this is the unique geodesic from identity to corresponding words.

\begin{enumerate}
\item
Suppose that there is a geodesic from $e$ to $\gamma(a^Nb^N)^k$ that
does not pass through $\gamma$.
Let $P_1$ and $P_3$ be the closest to $\gamma$ intersection points
of $\gamma(a^Nb^N)^k$ with this geodesic. Let $P_2=\gamma$.
Consider a triangle $P_1P_2P_3$ ($P_1P_2$ is some end subword$w$
of $\gamma$ and $P_2P_3$ is a beginning of $(a^Nb^N)^k$.
Let $w$ be the word that corresponds to the boundary of this triangle.
Note that $w$ is freely reduced and that $w$ (and any word
that we can get by symmetrization out of $w$ ) contains
at most 2 non-overlapping words which are $>3/6W$ and
at most one which is $>5/6W$.
This is a contradiction to Greendlinger lemma \cite{LynS} .
\item
First we see that by Greendlinger lemma that $b^k\gamma_2$ is geodesic
and that any geodesic from $e$ to $b^k\gamma_2$ passes through
$b^k$ (Similary to 1.).

Suppose that some geodesic from $e$ to $\gamma_1\gamma_2$ does not
pass through $\gamma_1$. Let $P_1$ and $P_3$ be the closest to $\gamma_1$
intersection points of this geodesic with $\gamma_1\gamma_2$.
Let $P_2=\gamma_1$. Consider the triangle $P_1P_2P_3$.
Note that this triangle is reduced (that is intersection points of its
sides are at $P_1, P_2, P_3$ only.
Note also that $P_1P_2$ contains $b^k$ as a subword.

Consider a reduced diagram of this triangle and its tiling.

{\bf First case.}
Suppose that its tiling contains at most 2 distinguished tiles.
Then we can apply part 3 of theorem 35 \cite{Stre} and see that
each tile that is not distinguished has 2 interior edges and 2 exterior
edges (see the picture in the proof of lemma 4).
 Note that there is a non-distinguished
 tile $B$ such that one of its exterior tiles lies within $b^k$.
Let $l$ be the length of the boundary of $B$. Note that each of
2 interior edges of $B$ have length at most $1/50l$.
The exterior tile that lies within $b^k$ has length at most  $1/50l$
and the other exterior tile has length $<1/2l$. We come to a contradiction.

{\bf Second case.}
Suppose that the tiling has 3 distinguished tiles.
Consider $P_2'=P_2a^{\pm 1}$. By Greendlinger lemma we see that
$P_1P_2'$ is geodesic.
Consider the triangle $P_1P_2'P_3$. It is a reduced triangle.
Consider a reduced diagram of this triangle and its tiling.
The tiling contains at most 2 distinguished tiles, since
none of $a^{\pm 1}b^{\pm 1}a^{\pm 1}$ is a piece.
In the same way as in the first case we come to a contradiction.
\end{enumerate}
{\bf Proof of the proposition.}
Let $M$ be the adjacency matrix of  Markov grammar of $H$.

Let $A_1, A_2,... A_f$ be irreducible blocks of $M$. Without loss
of generality we can assume that $M$ is
equal to
$$
\begin{array}{cccc}
A_1 & 0 & \ldots&0\\
* & A_2 & \ldots & 0\\
\vdots & \vdots &\ddots &\vdots \\
* & * & \ldots &A_f\\
\end{array}
$$
 The part of the Markov grammar that corresponds to $A_i$
will be called  $A_i$.
We say that $A_i$ is important if $v(A_i)=v(M)$.

We say that $A_i$ is good if we can generate $(a^Nb^N)^{100}$ inside
this block (starting at some state $a$ and ending at some state $b$ of $A_i$).

Consider the language generated by $A_i$ with the word  $(a^Nb^N)^{100}$
forbidden. Let
 $\tilde{v}(A_i)$ be the exponential growth rate of this language.

If $A_i$ is good and $v(A_i)>0$ then  $\tilde{v}(A_i)<v(A_i)$. In fact,
$A_i^{200N}$ is irreducible. Consider a matrix $A$ such that
$$
A(x,y)=A_i^{200N}(x,y)
$$
for $x\ne a$ or $y\ne b$ and
$$
A(a,b)=A_i^{200N}(a,b)-1.
$$
Since $A_i$ and $A$ are both non-negative
 we can deduce from lemma 2.3.2 of the chapter XIII \cite{Gant} that
$$
v(A)<v(A_i^{200N}).
$$
Hence
$$
\tilde{v}(A_i)^{200N}\le v(A)<v(A_i^{200N})=v(A_i)^{200N}.
$$
So we have shown that $\tilde{v}(A_i)<v(A_i)^{200N}$.

Note that $\tilde{v}(M)=\max \tilde{v}(A_i)$.
Hence it suffices to show that all important blocks are good.
Suppose there exist blocks that are important but not good.
Consider maximal $i$ such that $A_i$ is important but not good.
Let us prove that
\begin{enumerate}
\item
There exist $j>i$ such that $A_j$ is good.
\item
Let $j$ ($j>i$) be the maximal $j$ such that $A_j$ is good.
Then
$\tilde{v}(A_j)\ge v(A_i)$.
\end{enumerate}

Note that 1 and 2 are sufficient for contradiction, since 2 implies that
$$
v(A_j)>\tilde{v}(A_j)\ge v(A_i)
$$
and this means that $A_i$ is not important.

{\bf Proof of 1.}
Since $v(A_i)>0$ there exist two edges $O_1O_2$ and $O_2O_3$ inside
$A_i$ such that $\overrightarrow{O_1O_2}$
 is marked with $a$ or $a^{-1}$ and
$\overrightarrow{O_2O_3}$ is  marked with $b$ or $b^{-1}$.
Hence there exists a word $\gamma$ starting at the initial
state of $M$ and ending in $A_i$ with $\overrightarrow{O_1O_2}$,
$\overrightarrow{O_2O_3}$
(and hence ending with $a^{\pm 1}b^{\pm 1}$).
By first part of the lemma we see that
$\gamma(a^N b^N)^k$ is geodesic for any positive $k$
and that any geodesic from $e$ to $\gamma(a^N b^N)^k$ passes
through $\gamma$ and coincides from this point with
$\gamma(a^N b^N)^k$.
This means that we can generate $(a^N b^N)^k$ starting from $O_3$.
Then for some $j\ge i$ we can generate $(a^N b^N)^{100}$ inside $A_j$ and
so $A_j$ is good. (Note that we can never get from $A_k$ to $A_m$ if
$k>m$.)

{\bf Proof of 2.} Take maximal $j>i$ such that $A_j$ is good.
Then there exist $O_1, O_2,.., O_{s'}$ inside $A_j$ such that all the
edges
 $\overrightarrow{O_1O_2}, \overrightarrow{O_2O_3},...,
\overrightarrow{O_{s'}O_1}$
 are marked with $b$.
(Note, that $N$ is greater than the number of states in the grammar).

Then there exists a word $\gamma_1$ starting at the initial state
of $M$ and ending with $b^k$ in $O_1$ (for some $k>N$).
Consider some word $\gamma_2$ generating inside $A_i$ and starting with
$a^{\pm 1}b^{\pm 1}$.
Applying the second part of the lemma we get that we can generate $\gamma_2$
starting from $O_1$.

Note that for any word $\gamma_3$ generated by $A_i$ there exists a word
$\gamma_2$ generated by $A_i$ such that $\gamma_2$ begins with
$a^{\pm 1}b^{\pm 1}$, $\gamma_3$ is a subword of $\gamma_2$ and
$$
|\gamma_2|-|\gamma_3|\le s.
$$
Let $M'$ be the matrix corresponding to the union of blocks $A_j, A_{j+1},...
A_f$.
That is $M'$ is equal to
$$
\begin{array}{cccc}
A_j & 0 & \ldots&0\\
* & A_{j+1} & \ldots & 0\\
\vdots & \vdots &\ddots &\vdots \\
* & * & \ldots &A_f\\
\end{array}
$$

 We proved that
$$
\tilde{v}(M')\ge v(A_i).
$$
Hence there exists $j'\ge j$ such that
$$
\tilde{v}(A_{j'})\ge v(A_i).
$$
Note that $j'=j$ since otherwise $j'>j$ and $A_{j'}$ is important.
This is impossible.
In fact, if $A_{j'}$ is good, then $j$ is not maximal
among those for which   $A_j$ is good.
And if  $A_{j'}$ is not good, then $i$ is not maximal
among those for which $A_i$ is
important but not good.

\section{Proof of the Proposition 1}
Consider $E:\N \to \N$ such that
$$
E(i+1)>400E(i),
$$

$$
400/\sqrt{E(i+1)}< \alpha(200E(i))
$$
and
$$
E(i+1)>\gamma(\beta(200E(i)),2)
$$
Suppose that $J\ne M$ are subsets of $\N$.
Let
$$
J=\{j_1,j_2,...,j_k,... \},
$$
$$
M=\{m_1,m_2,...,m_k,... \},
$$
 $j_1<j_2<...<j_k<... $ and $m_1<m_2<...
<m_k<...$.

Suppose that $j_i=m_i$ for $i\le k$ and that
$j_{k+1}\ne m_{k+1}$.
Without loss of generality we can assume that $j_{k+1}>m_{k+1}$.
Let
$$
L=\{j_1, j_2,...,j_k^\}.
$$

 Let $v=v(G_L)$, $v_1=v(G_{J})$ and
$v_2=v(G_{M})$.
From the estimate from below and from lemma 1 we conclude that
$$
v_1\ge v-200/\sqrt{E(j_{k+1})}-200/\sqrt{E(j_{k+2})}
-200/\sqrt{E(j_{k+3})}-...
$$

$$
\ge v-400/\sqrt{E(j_{k+1})}>v-\alpha(200E(j_{k+1}-1))\ge
v-\alpha(200E(m_{k+1})).
$$
Let
$$
L'=\{m_1, m_2,...,m_k^{(2)}, m_{k+1}\}
$$
and let $v'=v_{G_{L'}}$.
From the estimate from above (proposition and definition of $\alpha$)
we get
$$
v'\le v-\alpha(200E(m_{k+1}))
$$
and hence
$$
v_2\le v'<v_1.
$$
So we have proved that $v_1 \ne v_2$.

\begin{rem}
In fact one can construct $E(i)$ described in proposition 1
explicitly. To do that one should estimate $\alpha(x)$ and
$\beta(x)$ for the case of small cancellation groups.
This can be done (for the example from the proof of the main theorem
of \cite{GyHa} one can
estimate the number of states of Markov grammar
 in terms of the hyperbolicity constant).
Ultimately one can see that
one can take
$$
E(i)=\underbrace{1000^{1000^{...1000}}}_{10i\quad {\rm times}}.
$$
We do not give the proof of this fact, since this
$E(i)$ grows very fast and because of that it seems to be of no practical use.
\end{rem}

I would like to thank Pierre de la Harpe for the encouragement and
many helpful comments on this paper. A part of this paper was written 
during the author's stay at the University of Geneva. The author 
gratefully aknowledges the support of the Swiss National Science Foundation.

\end{document}